\documentclass[12pt]{article}
\usepackage{amssymb}
\usepackage[margin=1.25 in]{geometry}
\usepackage{amsmath,epsfig}
\usepackage{times}
\usepackage[numbers]{natbib}
\usepackage{caption}
\usepackage{color}

\begin{document}

{\noindent{\large\bf{Correlation with Applications (1):}}

\vspace{3mm}

{\centering{\Large\bf{Measures of Correlation for Multiple
Variables}}

%}

\vspace{0.3 cm} {Jianji Wang, Nanning Zheng}

{\small{Institute of Artificial Intelligence and Robotics, Xi'an
Jiaotong University, Xi'an, China, 710049

wangjianji@mail.xjtu.edu.cn, nnzheng@mail.xjtu.edu.cn}}

}

{\center{\bf{Abstract}}

}

Multivariate correlation analysis plays an important role in
various fields such as statistics, economics, and big data
analytics. In this paper, we propose a pair of measures, the
unsigned correlation
coefficient (UCC) and the unsigned incorrelation coefficient (UIC), to
measure the strength of correlation and incorrelation (lack of correlation)
among multiple variables. The absolute value of Pearson's correlation coefficient is a special case of
UCC for two variables. Some important properties of UCC and UIC show that the proposed UCC and UIC are a pair of effective
measures for multivariate correlation. We also take the unsigned tri-variate correlation coefficient as an example to visually
display the effectiveness of the proposed UCC, and the geometrical explanation of UIC is also discussed. 
All the properties and the figures of UCC and UIC show that the proposed UCC and UIC are the general measures of correlation for multiple variables.

%We also take the UCC for three variables as an example to visually verify the effectiveness of the proposed UCC, and the geometrical explanation of UIC is also discussed.
%, and a new explanation of determinant is also made from the view of multivariate correlation.

\vspace{0.2 cm} {\bf{Key words}}: Unsigned correlation
coefficient (UCC), Unsigned incorrelation coefficient (UIC), Correlation
matrix, Strength of correlation.

\vspace{0.2 cm}

{\center\bf{1. \quad INTRODUCTION}

}

Correlation analysis is a statistical \,subject which studies
linear relationship and the ``strength" of linear relationship among
variables. It has been widely applied not only in statistics but
also in almost all fields of science \cite{Stigler}.

The research on quantitative method of correlation is one of the main research strategies for
correlation analysis, in which the strength of correlation
among variables is measured by a correlation measurement. Pearson's correlation coefficient is
a well-known bivariate correlation measurement.

For vectors ${\bf{a}}=(a_1,a_2,\cdots,a_n)^T$ and ${\bf{b}}=(b_1,b_2,\cdots,b_n)^T$, let $\mu_{\bf{a}}=\frac{1}{n}{\sum\limits_{i=1}^n{a_i}}$ and $\mu_{\bf{b}}=\frac{1}{n}{\sum\limits_{i=1}^n{b_i}}$ be the means of the elements in ${\bf{a}}$ and ${\bf{b}}$, respectively, let $\sigma^2_{\bf{a}}=\frac{1}{n-1}{\sum\limits_{i=1}^n}(a_i-\mu_{\bf{a}})^2$ and $\sigma^2_{\bf{b}}=\frac{1}{n-1}{\sum\limits_{i=1}^n}(b_i-\mu_{\bf{b}})^2$ be the variances of the elements in ${\bf{a}}$ and ${\bf{b}}$, respectively, and let $\sigma _{\bf{ab}}=\frac{1}{n-1}{\sum\limits_{i=1}^n}(a_i-\mu_{\bf{a}})(b_i-\mu_{\bf{b}})$ be the covariance between ${\bf{a}}$ and ${\bf{b}}$. If $\sigma _{\bf{a}}$ is not zero, we call ${\bf{a}}$ `a non-zero-variance vector'. Then for two non-zero-variance variables
${\bf{a}}$ and ${\bf{b}}$, Pearson's correlation coefficient between them is defined as
\begin{equation}
\rho_{{\bf{ab}}}  = \frac{{\sigma _{{\bf{ab}}} }}{{\sigma _{\bf{a}}
\sigma _{\bf{b}} }}= \cos \gamma.
\end{equation}
where $\gamma$ is the angle between the zero-mean variables of ${\bf{a}}$ and ${\bf{b}}$. In this paper, we call the variable ${\bf{a}}-\mu_{\bf{a}}{\bf{1}}$ the zero-mean variable of ${\bf{a}}$, where ${\bf{1}}$ is the vector with all ones. Pearson's
correlation coefficient $\rho _{{\bf{ab}}}\in [-1,1]$ and its absolute
value $|\rho _{{\bf{ab}}}|$ can be used to measure the strength of
correlation between ${\bf{a}}$ and ${\bf{b}}$. Let $\gamma^*$
be the angle between the zero-mean variables of ${\bf{a}}$ and ${\bf{b}}$ without considering
their directions, we have
\[
\gamma^* = \left\{ {\begin{array}{*{20}c}
   {\gamma ,\quad\qquad\ \gamma  \le 90^\circ}  \\
   {180^\circ  - \gamma ,\ \gamma  > 90^\circ}  \\
\end{array}} \right. ,
\]
and \[|\rho_{{\bf{ab}}}|  = \cos \gamma^*.\]

The linear relation of the variables ${\bf{a}}$ and ${\bf{b}}$ then
can be completely determined by the angle $\gamma^*$. If
$\gamma^*=0^\circ$, ${\bf{a}}$ and ${\bf{b}}$ are linear
dependent; If $\gamma^* = 90^\circ$, ${\bf{a}}$ and ${\bf{b}}$ are
perpendicular to each other and have the minimum correlation.

\begin{figure}[t]
\centering \centerline{\epsfig{figure=Fig1.pdf,width= 5.4 cm}}
\caption{The curve of the absolute value of Pearson's correlation
coefficient against the angle $\gamma$ between two zero-mean
variables.}\label{Fig1}
\end{figure}

It shows the curve of $|\rho_{{\bf{ab}}}|$ against the angle $\gamma$ within the intervals
$[0^\circ,180^\circ]$ in Figure~\ref{Fig1}, from which we can see that the absolute value of Pearson's correlation coefficient is a good
bivariate correlation measurement for two reasons. Firstly, the curve of $|\rho _{{\bf{ab}}}|$ against the angle $\gamma$ is always kept the opposite trend of change of the angle $\gamma^*$. Secondly, the strength of correlation measured by the absolute value of Pearson's correlation coefficient (the black curve in Figure~\ref{Fig1}) is close to the strength of correlation measured by the angle $\gamma^*$ (the blue dotted segments in Figure~\ref{Fig1}).

%The absolute value of Pearson's correlation coefficient is a good
%bivariate correlation measurement because the curve of $|\rho _{{\bf{ab}}}|$ against the angle $\gamma^*$
%is always kept the same trend of change with the line segments connecting the extreme-value
%points $|\rho_{{\bf{ab}}}|  = 1$ and $|\rho_{{\bf{ab}}}|=0$, and is always kept the opposite trend of the change of $\gamma^*$. Figure~\ref{Fig1} shows the curve of $|\rho
%_{{\bf{ab}}}|$ against the angle $\gamma$ within the intervals
%$[0^\circ,90^\circ]$ and $[90^\circ,180^\circ]$.

Pearson's correlation coefficient can only describe the correlation
relationship between two variables. Lots of
applications, such as data dimensionality reduction, subset selection, and sparse regression, need a correlation measurement to
measure the strength of correlation among multiple variables.
These applications usually select several variables from a number of
variables to make the selected variables have the minimum
correlation while maintaining enough information.

The development of information technology and big data analytics further
increases the importance of multivariate correlation
analysis. Unfortunately, there is no compact formulation to define and measure correlation for multiple variables. People have to estimate the strength of multivariate correlation by the indirect methods such as the partial correlation and the coefficient of determination. However, for most systems, such as physical, sociological, and economic systems, it is more important to discover associations among more than two variables.

In this paper, an unsigned correlation coefficient (UCC) and an unsigned
incorrelation coefficient (UIC) are proposed.
The proposed UCC and UIC can be used to measure the strength of
multivariate correlation and linear irrelevance. Many
properties of them are introduced in the paper. For example, both
the proposed UCC and UIC belong to an interval [0,1]; The sum of the
squared UCC and the squared UIC of a group of variables is always 1; The
value of UCC for multiple non-zero-variance variables
achieves the maximum value 1 if and only if these variables are
linear dependent, and achieves the minimum value 0 if and only if
these variables are perpendicular to each other; The value of UCC
for a group of variables is not less than the value of UCC for part of them; If
the number of variables is the same as their dimension, their UIC
then equals to the absolute value of the determinant of the square
matrix whose row or column vectors are the standardized vectors of
these variables.

We also show that the value of UIC for multivate variables is the volume of the parallelotope formed by the standardized vectors of them. Then we visually display that the strength of correlation of three variables measured by the unsigned tri-variate correlation coefficient is very close to the strength of correlation measured by a spatial angle $\eta$. 

All the properties and the figures of UCC and UIC show that the proposed UCC and UIC are the general measures of correlation for multiple variables.

\vspace{0.2 cm}

{\center\bf{2. \quad INNER PRODUCT-DETERMINANT EQUATION}

}

If two non-zero-variance variables ${\bf{a}}$ and ${\bf{b}}$ are linear dependent, we have $\rho^2_{{\bf{ab}}}= 1$. Garnett had proved that if three non-zero-variance variables ${\bf{a}}$,  ${\bf{b}}$, and ${\bf{c}}$  are linear dependent, then $\rho_{{\bf{ab}}}^2  + \rho_{{\bf{bc}}}^2  +
\rho_{{\bf{ac}}}^2 - 2\rho_{{\bf{ab}}} \rho_{{\bf{bc}}} \rho_{{\bf{ac}}}=1$ \cite{Garnett}.

For two variables ${\bf{a}}$ and ${\bf{b}}$, we use $|\rho_{{\bf{ab}}}|$ to measure the strength of correlation between them. According to the above discussion, we analyze the relation of the strength of correlation among ${\bf{a}}$,  ${\bf{b}}$, and ${\bf{c}}$  with the value of $r_{{\bf{abc}}}=\sqrt{\rho_{{\bf{ab}}}^2  + \rho_{{\bf{bc}}}^2  +
\rho_{{\bf{ac}}}^2 - 2\rho_{{\bf{ab}}} \rho_{{\bf{bc}}} \rho_{{\bf{ac}}}}$.

According to Garnett's conclusion, $r_{{\bf{abc}}}=1$ if and only if variables
${\bf{a}}$, ${\bf{b}}$, and ${\bf{c}}$ on the same plane. Moreover, it is obviously that $r_{{\bf{abc}}}=0$ if and only if variables
${\bf{a}}$, ${\bf{b}}$, and ${\bf{c}}$ are perpendicular to each
other. The two important properties of $r_{{\bf{abc}}}$ show that $r_{{\bf{abc}}}$ may be a proper measurement to measure the strength of correlation among three variables.

We denote by ${\bf{a}}'$ the standardized vector of
the non-zero-variance vector ${\bf{a}}$:
\[{\bf{a}}' = \frac{{\bf{a}}-{\mu_{\bf{a}}}{\bf{1}}}{\left\|{\bf{a}}-{\mu_{\bf{a}}}{\bf{1}}\right\|}.\]

If ${\bf{a}}{\rm{'}}$ and ${\bf{b}}{\rm{'}}$ are the standardized
vectors of ${\bf{a}}$ and ${\bf{b}}$, respectively, Pearson's
correlation coefficient between ${\bf{a}}$ and ${\bf{b}}$ is then
the inner product between ${\bf{a}}{\rm{'}}$ and ${\bf{b}}{\rm{'}}$, and $r_{{\bf{abc}}}$ can also be expressed by the inner product of ${\bf{a}}{\rm{'}}$ and ${\bf{b}}{\rm{'}}$, the inner product of ${\bf{a}}{\rm{'}}$ and ${\bf{c}}{\rm{'}}$, and the inner product of ${\bf{b}}{\rm{'}}$ and ${\bf{c}}{\rm{'}}$.

Here we focus on an identical relation between inner
product and a determinant group, which we call the inner
product-determinant equation (IPD equation).

For $n$-dimensional variables ${\bf{a}}_1$, ${\bf{a}}_2$, $\cdots$,
and ${\bf{a}}_m$, ${\bf{a}}_i = (a_{i1} ,a_{i2} , \cdots$, $a_{in}
)^T$, $i = 1, 2, \cdots$, $m$, $m \le n$, we denote
\[
[{\bf{a}}(m)|j_1 ,j_2 , \cdots ,j_m ] = \left[
{\begin{array}{*{20}c}
   {a_{1j_1 } } & {a_{1j_2 } } &  \cdots  & {a_{1j_m } }  \\
   {a_{2j_1 } } & {a_{2j_2 } } &  \cdots  & {a_{2j_m } }  \\
    \vdots  &  \vdots  &  \ddots  &  \vdots   \\
   {a_{mj_1 } } & {a_{mj_2 } } &  \cdots  & {a_{mj_m } }  \\
\end{array}} \right],
\]
where $1 \le j_1  < j_2  <  \cdots  < j_m  \le n$.

According to Cauchy-Binet formula \cite{Zeng}, we have the following
lemma:

\noindent{\bf{Lemma 1}} \ For $n$-dimensional variables ${\bf{a}}_1$,
${\bf{a}}_2$, $\cdots$, and ${\bf{a}}_m$, ${\bf{a}}_i = (a_{i1}
,a_{i2} , \cdots, a_{in} )^T$, $i\in\{1,2,\cdots,m\}$, $m\leq n$, if
${\bf{M}}$ is the inner product matrix of these variables,
\[{\bf{M}} = \left[ {\begin{array}{*{20}c}
   {<\!{\bf{a}}_1 ,{\bf{a}}_1 \!>} & {<\!{\bf{a}}_1 ,{\bf{a}}_2 \!>} &  \cdots  & {<\!{\bf{a}}_1 ,{\bf{a}}_m \!>}  \\
   {<\!{\bf{a}}_2 ,{\bf{a}}_1 \!>} & {<\!{\bf{a}}_2 ,{\bf{a}}_2 \!>} &  \cdots  & {<\!{\bf{a}}_2 ,{\bf{a}}_m \!>}  \\
    \vdots  &  \vdots  &  \ddots  &  \vdots   \\
   {<\!{\bf{a}}_m ,{\bf{a}}_1 \!>} & {<\!{\bf{a}}_m ,{\bf{a}}_2 \!>} &  \cdots  & {<\!{\bf{a}}_m ,{\bf{a}}_m \!>}  \\
\end{array}} \right] ,
\]
then we have
\begin{equation}
\det ({\bf{M}}) = \sum\limits_{j_1  < j_2  <  \cdots < j_m } {(\det
[{\bf{a}}(m)|j_1 ,j_2 , \cdots ,j_m ])^2 },
\end{equation}
where $j_1 ,j_2 , \cdots ,j_m \in \{1,2,\cdots,n\}$.

Proof: According to Cauchy-Binet formula, for the matrixes
${\bf{A}}=[{\bf{a}}_1,{\bf{a}}_2, \cdots,{\bf{a}}_m]$ and
${\bf{A}}^T$, we have
\[
\det ({\bf{A}}^T{\bf{A}}) = \sum\limits_{j_1  < j_2  <  \cdots  <
j_m } {(\det [{\bf{a}}(m)|j_1  < j_2  <  \cdots  < j_m ])^2 } .
\]

Because ${\bf{A}}^T {\bf{A}}={\bf{M}}$, this lemma is true.\hfill$\blacksquare$

We denote by $\Gamma$ the circular inner product
\begin{equation}
\begin{array}{l}
\Gamma\!\!<\!{\bf{a}}_{k_1} ,{\bf{a}}_{k_2}, \cdots ,{\bf{a}}_{k_p}\! >
= <\!{\bf{a}}_{k_1},{\bf{a}}_{k_2}\!><\!{\bf{a}}_{k_2},{\bf{a}}_{k_3}\!>
\cdots <\!{\bf{a}}_{k_{p-1}} ,{\bf{a}}_{k_p}
\!><\!{\bf{a}}_{k_p},{\bf{a}}_{k_1} \!>.
\end{array}
\end{equation}

For a group of variables, all the permutations of them which can generate
the same inner product or circular inner product are regarded as the
same inner-product-permutation, then the inner product-determinant
equation can be rewritten as following:

\noindent{\bf{Inner product-Determinant Equation}} \ For
$n$-dimensional variables ${\bf{a}}_1$, ${\bf{a}}_2$, $\cdots$, and
${\bf{a}}_m$, $m\leq n$, the inner-product-determinant equation is
\begin{equation}
\begin{array}{l}
\sum\limits_\pi  {2^{\left| {\pi _3 } \right|} ( - 1)^{m - \left|
{\pi} \right|} }{\prod\limits_{\pi _1 } \left\| {{\bf{a}}_s } \right\|^2\prod\limits_{\pi _2 } {<\!{\bf{a}}_i ,{\bf{a}}_j\!>^2 } }
\prod\limits_{\pi _3 } {\Gamma\!\!<\!{\bf{a}}_{k_1} ,{\bf{a}}_{k_2} ,
\cdots ,{\bf{a}}_{k_p }\!>}  \\
= \sum\limits_{j_1  < j_2  <  \cdots < j_m } {(\det [{\bf{a}}(m)|j_1
,j_2 , \cdots ,j_m ])^2 },
\end{array}
\end{equation}
where $\pi$ runs through the list of all partitions of
$\{{\bf{a}}_1, {\bf{a}}_2, \cdots$, ${\bf{a}}_m\}$ with only
consider the inner-product-permutation, and the subsets in each
partition $\pi$ are divided into three classes: If one subset only
contains one variable, then this subset belongs to the first class
$\pi_1$, and the self-inner product of the variable appears in the
formula of the partition; If one subset contains two different variables, then
this subset belongs to the second class $\pi_2$, and the square of
mutual inner product between the two variables appears in the
formula of the partition; The others belong to the third class
$\pi_3$, and the circular inner product for each subset in $\pi_3$
appears in the formula of the partition; The number of subsets in a
partition $\pi$ and $\pi_3$ are $\left| {\pi} \right|$ and $\left|
{\pi _3 } \right|$, respectively.

Proof: In fact, each partition $\pi$ in the inner
product-determinant equation is corresponding to one item in the
expansion of the determinant of the inner product matrix. By
considering the inversion number and the symmetry of inner product
matrix, the above equation can be easily obtained.\hfill$\blacksquare$

If the number of subsets in each $\pi_1$ and $\pi_2$ are $\left|
{\pi_1 } \right|$ and $\left| {\pi _2 } \right|$, respectively, then
$\left| {\pi_1 } \right|+\left| {\pi_2 } \right|+\left| {\pi_3 }
\right|=\left| {\pi } \right|$. Moreover, the subsets in $\pi_3$ are
not the standard sets because the inner-product-permutation is
involved.

The inner product-determinant equation has a set-partition-based
form, which is similar to the joint cumulant equation
\cite{cumulant, cumulants}.

For three $n$-dimensional variables ${\bf{a}}$, ${\bf{b}}$, and
${\bf{c}}$, $n\ge 3$, there are 6 cases of inner products in total,
and they are $<\!{\bf{a}},{\bf{a}}\!>$, $<\!{\bf{b}},{\bf{b}}\!>$,
$<\!{\bf{c}},{\bf{c}}\!>$, $<\!{\bf{a}},{\bf{b}}\!>$, $<\!{\bf{a}},{\bf{c}}\!>$,
and $<\!{\bf{b}},{\bf{c}}\!>$. From Lemma 1, the IPD equation for three
variables can be obtained as following:

\noindent{\bf{Corollary 1}} \ For three $n$-dimensional variables
${\bf{a}}$, ${\bf{b}}$, and ${\bf{c}}$, $n\geq 3$, the inner
product-determinant equation is as following:
\begin{equation}
\begin{array}{l}
\left\| {\bf{a}} \right\|^2 \left\| {\bf{b}} \right\|^2 \left\| {\bf{c}}
\right\|^2   - \left\|
{\bf{a}} \right\|^2 \!<\!{\bf{b}},{\bf{c}}\!>^2\!  - \left\| {\bf{b}}
\right\|^2 \!<\!{\bf{a}},{\bf{c}}\!>^2\! - \left\| {\bf{c}} \right\|^2\!
<\!{\bf{a}},{\bf{b}}\!>^2\! \\+
2<\!{\bf{a}},{\bf{b}}\!><\!{\bf{b}},{\bf{c}}\!><\!{\bf{a}},{\bf{c}}\!>\ = \sum\limits_{i < j < k} {<\!\det
[{\bf{a}},{\bf{b}},{\bf{c}}|i,j,k]\!>^2 },
\end{array}
\end{equation}
where $i,j,k \in \{1,2, \cdots, n\}$.

Similarly, IPD equation for two variables is listed below.

\noindent{\bf{Corollary 2}} \ For two $n$-dimensional variables
${\bf{a}}$ and ${\bf{b}}$, $n\geq 2$, the inner product-determinant
equation is
\begin{equation}
\left\| {\bf{a}} \right\|^2 \left\| {\bf{b}} \right\|^2  -
<\!{\bf{a}},{\bf{b}}\!>^2  = \sum\limits_{i < j} {(\det
[{\bf{a}},{\bf{b}}|i,j])^2 },
\end{equation}
where $i,j \in \{1,2, \cdots, n\}$.

\vspace{2mm}

{\center\bf{3. \quad CORRELATION MEASURES FOR MULTIPLE VARIABLES}

}

According to Corollary 1, if ${\bf{a}}{\rm{'}}$, ${\bf{b}}{\rm{'}}$, and
${\bf{c}}{\rm{'}}$ are the standardized vectors of ${\bf{a}}$,
${\bf{b}}$, and ${\bf{c}}$, respectively, we have
\begin{equation}
r_{{\bf{abc}}}^2  = 1 - \sum\limits_{i  < j  < k } {(\det
[{\bf{a}}',{\bf{b}}',{\bf{c}}'|i ,j ,k ])^2 } .
\end{equation}

Similarly, the square of Pearson's correlation coefficient between
${\bf{a}}$ and ${\bf{b}}$ can be obtained from Corollary 2:
\begin{equation}
\rho_{{\bf{ab}}}^2  = 1 - \sum\limits_{i  < j } {(\det
[{\bf{a}}',{\bf{b}}'|i ,j ])^2 } .
\end{equation}

Inspired by the above formulas of $r_{{\bf{abc}}}^2$ and
$\rho_{{\bf{ab}}}^2$, we define the multivariate correlation
measurement as following:

\noindent{\bf{Definition 1}} \ For $n$-dimensional non-zero-variance
variables ${\bf{a}}_1$,  ${\bf{a}}_2$, $\cdots$, ${\bf{a}}_m$, $2\le
m \le n$, if ${\bf{a}}'_i=(a'_{i1},a'_{i2},\cdots,a'_{in})^T$ is the standardized vector of
${\bf{a}}_i=(a_{i1},a_{i2},\cdots,a_{in})^T$, $i\in\{1,2,\cdots,m\}$, then the unsigned correlation
coefficient $r_{{\bf{a}}_1 {\bf{a}}_2 \cdots {\bf{a}}_m }$ among
${\bf{a}}_1$, ${\bf{a}}_2$, $\cdots$, and ${\bf{a}}_m$ is defined as
\begin{equation}
\begin{array}{l}
 r_{{\bf{a}}_1 {\bf{a}}_2  \cdots {\bf{a}}_m }^2
 = 1 - \sum\limits_{j_1  < j_2  <  \cdots  < j_m } {(\det [{\bf{a}}'(m)|j_1 ,j_2 , \cdots ,j_m ])^2 }\\
 \end{array},
\end{equation}
where
\[
[{\bf{a}}'(m)|j_1 ,j_2 , \cdots ,j_m ] = \left[
{\begin{array}{*{20}c}
   {a'_{1j_1 } } & {a'_{1j_2 } } &  \cdots  & {a'_{1j_m } }  \\
   {a'_{2j_1 } } & {a'_{2j_2 } } &  \cdots  & {a'_{2j_m } }  \\
    \vdots  &  \vdots  &  \ddots  &  \vdots   \\
   {a'_{mj_1 } } & {a'_{mj_2 } } &  \cdots  & {a'_{mj_m } }  \\
\end{array}} \right],
\]
$j_1, j_2, \cdots, j_m \in \{1,2,\cdots,n\}$ and $1 \le j_1  < j_2  <  \cdots  < j_m  \le n$.

The sign of mutual direction for two variables can be judged by
whether the angle between the zero-mean variables of them is larger than $90^\circ$. However, there is
no the mutual direction for multiple variables. Therefore, the
correlation measurement for multiple variables in this paper is defined as an unsigned value, the rationality of which is also discussed in Section 3 in this paper.

According to Definition 1, for two non-zero-variance
variables ${\bf{a}}$ and ${\bf{b}}$, $r_{{\bf{ab}}}=|\rho_{{\bf{ab}}}|$. Then we can see that the sum of the squares of the
determinant group is a coupling part of the proposed
UCC. We define the coupling part as the square of the unsigned
incorrelation coefficient (UIC), which can be used to measure linear
irrelevance among variables:

\noindent{\bf{Definition 2}} \ For $n$-dimensional non-zero-variance
variables ${\bf{a}}_1$, ${\bf{a}}_2$, $\cdots$, ${\bf{a}}_m$, $2\le
m \le n$, if ${\bf{a}}'_i$ is the standardized vector of
${\bf{a}}_i$, $i\in\{1,2,\cdots,m\}$, the unsigned incorrelation
coefficient (UIC) $\omega_{{\bf{a}}_1 {\bf{a}}_2 \cdots {\bf{a}}_m
}$ among ${\bf{a}}_1$, ${\bf{a}}_2$, $\cdots$, and ${\bf{a}}_m$ is
defined as
\begin{equation}
\begin{array}{l}
\omega_{{\bf{a}}_1 {\bf{a}}_2 \cdots {\bf{a}}_m }^2 =
\sum\limits_{j_1  < j_2  < \cdots  < j_m } {(\det [{\bf{a}}'(m)|j_1
,j_2 , \cdots ,j_m ])^2 },
\end{array}
\end{equation}
where $j_1, j_2, \cdots, j_m \in \{1,2,\cdots,n\}$ and $1 \le j_1  < j_2  <  \cdots  < j_m  \le n$.

A lemma exists for UIC as following:

\noindent{\bf{Lemma 2}} \ For $n$-dimensional variables ${\bf{a}}_1$,
${\bf{a}}_2$, $\cdots$, ${\bf{a}}_m$, ${\bf{a}}_{m+1}$, $2\le {m+1}
\le n$, if ${\bf{a}}'_i$ is the standardized vector of ${\bf{a}}_i$,
$i = 1,2, \cdots ,m, m+1$, we have
\begin{equation}
\begin{array}{l}
\omega _{{\bf{a}}_1 {\bf{a}}_2  \cdots {\bf{a}}_m }^2 - \omega
_{{\bf{a}}_1 {\bf{a}}_2  \cdots {\bf{a}}_m {\bf{a}}_{m + 1 } }^2 \\=
\sum\limits_{J_{m - 1} } {(\sum\limits_{p \notin J_{m - 1}} {( -
1)^{g(p:J_{m - 1} )} a'_{m + 1,p} \det [{\bf{a}}'(m)|p,J_{m - 1} ]}
)^2 }
\end{array},
\end{equation}
where $J_{m - 1}=\{j_1,j_2,\cdots,j_{m-1}\}$,
$j_1$$<$$j_2$$<$$\cdots$$<$$j_{m-1}$, $j_1$, $j_2$, $\cdots$,
$j_{m-1}$ $\in \{1,2, \cdots, n\}$, and $g(p: S)$ is the number of
elements which are larger than $p$ in the set $S$.

Proof: Let $J_m=\{j^\dagger_1,j^\dagger_2,\cdots,j^\dagger_m\}$,
$j^\dagger_1$$<$$j^\dagger_2$$<$$\cdots$$<$$j^\dagger_{m}$, $j^\dagger_1$, $j^\dagger_2$, $\cdots$, $j^\dagger_{m}$
$\in \{1$, $2, \cdots, n\}$, and $J_{m+1}=\{j^\ddagger_1,j^\ddagger_2,\cdots$, $j^\ddagger_m,
j^\ddagger_{m+1}\}$, $j^\ddagger_1$$<$$j^\ddagger_2$$<$$\cdots$$<$$j^\ddagger_{m}$$<$$j^\ddagger_{m+1}$, $j^\ddagger_1$,
$j^\ddagger_2$, $\cdots$, $j^\ddagger_{m}$, $j^\ddagger_{m+1}$ $\in \{1$, $2, \cdots, n\}$.
\[
\begin{array}{l}
 \omega _{{\bf{a}}_1 {\bf{a}}_2  \cdots {\bf{a}}_m {\bf{a}}_{m + 1} }^2  = \sum\limits_{J_{m+1}} {(\det [{\bf{a}}'(m + 1)|J_{m+1}])^2 }  \\
  = \sum\limits_{J_{m+1}} {\{ \sum\limits_{p = j_1 }^{j_{m + 1} } {({a'}_{m + 1,p}^2 \det ^2 [{\bf{a}}'(m)|J_{m+1} \backslash p])} }  \\
  + \sum\limits_{\scriptstyle l,p \in J_{m+1} \hfill \atop \scriptstyle \quad\, l \ne p \hfill} \begin{array}{l}
 [2( - 1)^{g(l:J_{m+1})}(a'_{m + 1,l} \det [{\bf{a}}'(m)|J_{m+1} \backslash l]) \\
 ( - 1)^{g(p:J_{m+1})}(a'_{m + 1,p} \det [{\bf{a}}'(m)|J_{m+1} \backslash p])] \\
 \end{array} \}.  \\
 \end{array}
\]

The first part can be rewritten as
\[
\begin{array}{l}
 \sum\limits_{J_{m+1}} {\sum\limits_{p = j_1 }^{j_{m + 1} } {({a'}_{m + 1,p}^2 \det ^2 [{\bf{a}}'(m)|J_{m+1} \backslash p])}}
 = \sum\limits_p {{a'}_{m + 1,p}^2 } \sum\limits_{\scriptstyle \ \, J_m  \hfill \atop
  \scriptstyle p \notin J_m  \hfill} {\det ^2 [{\bf{a}}'(m)|J_m ]}  \\
  = \sum\limits_p {{a'}_{m + 1,p}^2 } (\omega _{{\bf{a}}_1 {\bf{a}}_2  \cdots {\bf{a}}_m}^2 - \sum\limits_{\scriptstyle \ \, J_m  \hfill \atop
  \scriptstyle p \in J_m \hfill} {\det ^2 [{\bf{a}}'(m)|J_m ])}  \\
  = \left\| {{\bf{a}}'_{m + 1} } \right\|^2 \omega _{{\bf{a}}_1 {\bf{a}}_2  \cdots {\bf{a}}_m}^2  - \sum\limits_p {{a'}_{m + 1,p}^2 } \sum\limits_{\scriptstyle \ \, J_m  \hfill \atop
  \scriptstyle p \in J_m \hfill} {\det ^2 [{\bf{a}}'(m)|J_m ]}  \\
  = \omega _{{\bf{a}}_1 {\bf{a}}_2  \cdots {\bf{a}}_m}^2  - \sum\limits_p {{a'}_{m + 1,p}^2 } \sum\limits_{\scriptstyle \ J_{m - 1}  \hfill \atop
  \scriptstyle p \notin J_{m - 1} \hfill} {\det ^2 [{\bf{a}}'(m)|p,J_{m - 1} ]},  \\
 \end{array}
\]
and the second part can be rewritten as
\[
\begin{array}{l}
 \sum\limits_{J_{m + 1} } {\sum\limits_{\scriptstyle l,p \in J_{m + 1} \hfill \atop
  \scriptstyle \quad\ l \ne p \hfill} \begin{array}{l}
 \{ 2( - 1)^{g(l:J_{m + 1}) + g(p:J_{m + 1})} a'_{m + 1,l} a'_{m + 1,p}  \\
 \det [{\bf{a}}'(m)|J_{m + 1} \backslash l] \det [{\bf{a}}'(m)|J_{m + 1} \backslash p]\}  \\
 \end{array} }  \\
  = 2\sum\limits_{l \ne p} {a'_{m + 1,l} a'_{m + 1,p} \sum\limits_{\scriptstyle \quad J_{m + 1}  \hfill \atop
  \scriptstyle \ l,p \in J_{m + 1} \hfill} \begin{array}{l}
 \{( - 1)^{g(l:J_{m + 1}) + g(p:J_{m + 1})}  \\
 \det [{\bf{a}}'(m)|J_{m + 1} \backslash l]
 \det [{\bf{a}}'(m)|J_{m + 1} \backslash p]\} \\
 \end{array} }  \\
  = 2\sum\limits_{l \ne p} {a'_{m + 1,l} a'_{m + 1,p} \sum\limits_{\scriptstyle \quad J_{m - 1}  \hfill \atop
  \scriptstyle l,p \notin J_{m - 1} \hfill} \begin{array}{l}
 \{( - 1)^{g(l:J_{m - 1} ) + g(p:J_{m - 1} ) + 1}  \\
 \det [{\bf{a}}'(m)|p,J_{m - 1} ]
 \det [{\bf{a}}'(m)|l,J_{m - 1} ]\} \\
 \end{array} }  \\
  =  - 2\sum\limits_{J_{m - 1} } {\sum\limits_{
  {\scriptstyle \quad\, l \ne p \hfill \atop
  \scriptstyle l,p \notin J_{m - 1} \hfill}} \begin{array}{l}
 \{( - 1)^{g(p:J_{m - 1} )}
 a'_{m + 1,p} \det [{\bf{a}}'(m)|p,J_{m - 1} ] \\
 ( - 1)^{g(l:J_{m - 1} )} a'_{m + 1,l} \det [{\bf{a}}'(m)|l,J_{m - 1} ]\} \\
 \end{array} }.  \\
 \end{array}
\]

Hence, we have
\[
\omega _{{\bf{a}}_1 {\bf{a}}_2  \cdots {\bf{a}}_m {\bf{a}}_{m + 1}
}^2  = \omega _{{\bf{a}}_1 {\bf{a}}_2  \cdots {\bf{a}}_m {\bf{a}}_m
}^2  - \phi ({\bf{a}}'(m + 1)),
\]
and $\phi ({\bf{a}}'(m + 1))$ can be expressed as
\[
\begin{array}{l}
 \phi ({\bf{a}}'(m + 1)) = \sum\limits_{J_{m - 1} } {\sum\limits_{p \notin J_{m - 1} } \begin{array}{l}
 {a'}_{m + 1,p}^2 \det ^2 [{\bf{a}}'(m)|p,J_{m - 1} ] \\
 \end{array} }  \\
  + 2\sum\limits_{J_{m - 1} } {\sum\limits_{
  {\scriptstyle \quad\, l \ne p \hfill \atop
  \scriptstyle l,p \notin J_{m - 1} \hfill}} \begin{array}{l}
 \{ ( - 1)^{g(p:J_{m - 1} )} a'_{m + 1,p} \det [{\bf{a}}'(m)|p,J_{m - 1} ] \\
  \cdot ( - 1)^{g(l:J_{m - 1} )} a'_{m + 1,l} \det [{\bf{a}}'(m)|l,J_{m - 1} ]\}  \\
 \end{array} }  \\
  = \sum\limits_{J_{m - 1} } {(\sum\limits_{p \notin J_{m - 1} } \begin{array}{l}
 ( - 1)^{g(p:J_{m - 1} )} a'_{m + 1,p}
 \det [{\bf{a}}'(m)|p,J_{m - 1} ] )^2. \\
 \end{array} }  \\
 \end{array}
\]\hfill$\blacksquare$

The proposed unsigned correlation coefficient and unsigned
incorrelation coefficient have some important properties,
several of which are discussed below.

{\bf{Property 3.1}} \ $r_{{\bf{a}}_1 {\bf{a}}_2 \cdots {\bf{a}}_m }$
and $\omega_{{\bf{a}}_1 {\bf{a}}_2 \cdots {\bf{a}}_m }$ are both the
symmetric functions of ${\bf{a}}_1$, ${\bf{a}}_2$, $\cdots$,
${\bf{a}}_m$.

{\bf{Property 3.2}} \ If $r_{{\bf{a}}_1 {\bf{a}}_2 \cdots {\bf{a}}_m
}$ and $\omega_{{\bf{a}}_1 {\bf{a}}_2 \cdots {\bf{a}}_m }$ are
UCC and UIC for non-zero-variance variables
${\bf{a}}_1$, ${\bf{a}}_2$, $\cdots$, ${\bf{a}}_m$, respectively,
and $r_{{\bf{a}}_1 {\bf{a}}_2 \cdots {\bf{a}}_{m-1} }$ and
$\omega_{{\bf{a}}_1 {\bf{a}}_2 \cdots {\bf{a}}_{m-1} }$ are
UCC and UIC for ${\bf{a}}_1$, ${\bf{a}}_2$, $\cdots$,
${\bf{a}}_{m-1}$, respectively, then
\[
\begin{array}{l}
r_{{\bf{a}}_1{\bf{a}}_2\cdots{\bf{a}}_m}^2 \ge
r_{{\bf{a}}_1{\bf{a}}_2\cdots{\bf{a}}_{m-1}}^2 \\
\omega_{{\bf{a}}_1{\bf{a}}_2\cdots{\bf{a}}_m}^2 \le
\omega_{{\bf{a}}_1{\bf{a}}_2\cdots{\bf{a}}_{m-1}}^2 \\
\end{array}
\]

Property 3.2 can be directly obtained from Lemma 2. It shows that
the value of UCC for some variables is not less than the value of UCC for
part of the variables.

{\bf{Property 3.3}}
\[
\begin{array}{l}
 0 \le r_{{\bf{a}}_1 {\bf{a}}_2  \cdots {\bf{a}}_m }^2  \le 1 \\
 0 \le \omega _{{\bf{a}}_1 {\bf{a}}_2  \cdots {\bf{a}}_m }^2  \le 1 \\
 \end{array}
\]

Proof: According to Definition 2 we have $\omega _{{\bf{a}}_1
{\bf{a}}_2 \cdots {\bf{a}}_m }^2\geq 0$. From Property 3.2, $\omega
_{{\bf{a}}_1 {\bf{a}}_2 \cdots {\bf{a}}_m }^2$ $\le$ $\cdots \le
\omega _{{\bf{a}}_i {\bf{a}}_j}^2 = 1 - \rho _{{\bf{a}}_i {\bf{a}}_j}^2
\le 1$, $i$$\neq$$j,\ i,j\in\{1,2,\cdots,m\}$. Because $r
_{{\bf{a}}_1 {\bf{a}}_2 \cdots {\bf{a}}_m }^2+\omega _{{\bf{a}}_1
{\bf{a}}_2 \cdots {\bf{a}}_m }^2=1$, we have $0 \le r_{{\bf{a}}_1
{\bf{a}}_2 \cdots {\bf{a}}_m }^2 \le 1$.\hfill$\blacksquare$

{\bf{Property 3.4}} \ $r_{{\bf{a}}_1 {\bf{a}}_2 \cdots {\bf{a}}_m
}^2=1 $ if and only if variables ${\bf{a}}_1, {\bf{a}}_2, \cdots$,
and ${\bf{a}}_m$ are linear dependent.

Proof:

Sufficiency: If ${\bf{a}}_1, {\bf{a}}_2, \cdots$, and ${\bf{a}}_m$
are linear dependent, $\det [{\bf{a}}'(m)|j_1 ,j_2$, $\cdots ,j_m
]=0$ for all cases of $j_1<j_2<\cdots<j_m$, $j_1 ,j_2 , \cdots
,j_m\in\{1,2,\cdots,n\}$.

Necessity: We denote by $ {\mathop {\bf{a}}\limits^ \wedge}_j =
{(a_{1j} ,a_{2j} , \cdots ,a_{mj} )^T }$ the $j$th column vector of
the matrix [${\bf{a}}(m)|1,2$, $\cdots$, $n$], $j=1,2,\cdots,n$.
Presume that $j_1,j_2,\cdots,j_m$ exist to make rank$\{ {\mathop
{\bf{a}}\limits^ \wedge}_{j_1},{\mathop {\bf{a}}\limits^
\wedge}_{j_2},\cdots,{\mathop {\bf{a}}\limits^ \wedge}_{j_m}\}=m$,
then $\det [{\bf{a}}(m)|j_1 ,j_2 , \cdots ,j_m ]\ne0$ and
$r_{{\bf{a}}_1 {\bf{a}}_2 \cdots {\bf{a}}_m }^2<1 $. Hence, the
presume is incorrect. We have rank$\{{\mathop {\bf{a}}\limits^
\wedge}_1,{\mathop {\bf{a}}\limits^ \wedge}_2,\cdots,{\mathop
{\bf{a}}\limits^ \wedge}_n\}$ $\le m-1$ if $r_{{\bf{a}}_1 {\bf{a}}_2
\cdots {\bf{a}}_m }^2=1 $. Because the row rank equals to the column
rank of the same matrix, we have
\[
\text{rank}\{{\bf{a}}_1, {\bf{a_2}}, \cdots,
{\bf{a}}_m\}=\text{rank}\{{\mathop {\bf{a}}\limits^
\wedge}_1,{\mathop {\bf{a}}\limits^ \wedge}_2,\cdots,{\mathop
{\bf{a}}\limits^ \wedge}_n\}\le m-1 .
\]\hfill$\blacksquare$

For four variables ${\bf{a}}_1, {\bf{a}}_2, {\bf{a}}_3$, and
${\bf{a}}_4$, let $\rho _{{\bf{a}}_i {\bf{a}}_j}$ be Pearson's
correlation coefficient between the variables ${\bf{a}}_i$ and
${\bf{a}}_j$, $i,j\in\{1,2,3,4\}$. Then ${\bf{a}}_1, {\bf{a}}_2,
{\bf{a}}_3$, and ${\bf{a}}_4$ are linear dependent if and only if
\begin{equation}
\begin{array}{l}
\rho _{{\bf{a}}_1 {\bf{a}}_2}^2+\rho_{{\bf{a}}_1 {\bf{a}}_3}^2+\rho
_{{\bf{a}}_1{\bf{a}}_4}^2+\rho _{{\bf{a}}_2 {\bf{a}}_3}^2+\rho_{{\bf{a}}_2
{\bf{a}}_4}^2+\rho _{{\bf{a}}_3{\bf{a}}_4}^2 -(\rho _{{\bf{a}}_1
{\bf{a}}_2}^2\rho _{{\bf{a}}_3{\bf{a}}_4}^2+\rho_{{\bf{a}}_1
{\bf{a}}_3}^2\rho_{{\bf{a}}_2 {\bf{a}}_4}^2+\rho
_{{\bf{a}}_1{\bf{a}}_4}^2\rho _{{\bf{a}}_2
{\bf{a}}_3}^2)\\-2(\rho_{{\bf{a}}_1 {\bf{a}}_3}\rho
_{{\bf{a}}_1{\bf{a}}_4}\rho _{{\bf{a}}_3{\bf{a}}_4}+\rho _{{\bf{a}}_1
{\bf{a}}_2}\rho _{{\bf{a}}_1{\bf{a}}_4}\rho_{{\bf{a}}_2 {\bf{a}}_4}+\rho
_{{\bf{a}}_1 {\bf{a}}_2}\rho_{{\bf{a}}_1 {\bf{a}}_3}\rho _{{\bf{a}}_2
{\bf{a}}_3}+\rho_{{\bf{a}}_2 {\bf{a}}_3}\rho_{{\bf{a}}_2 {\bf{a}}_4}\rho
_{{\bf{a}}_3{\bf{a}}_4})\\+2(\rho_{{\bf{a}}_1 {\bf{a}}_3}\rho
_{{\bf{a}}_1{\bf{a}}_4}\rho _{{\bf{a}}_2 {\bf{a}}_3}\rho_{{\bf{a}}_2
{\bf{a}}_4}+\rho _{{\bf{a}}_1 {\bf{a}}_2}\rho _{{\bf{a}}_1{\bf{a}}_4}\rho
_{{\bf{a}}_2 {\bf{a}}_3}\rho _{{\bf{a}}_3{\bf{a}}_4}+\rho _{{\bf{a}}_1
{\bf{a}}_2}\rho_{{\bf{a}}_1 {\bf{a}}_3}\rho_{{\bf{a}}_2 {\bf{a}}_4}\rho
_{{\bf{a}}_3{\bf{a}}_4}) = 1
\end{array},
\end{equation}
which is also kept the same with Garnett's results \citep{Garnett}.

{\bf{Property 3.5}} \ $r_{{\bf{a}}_1 {\bf{a}}_2 \cdots {\bf{a}}_m
}^2=0 $ if and only if variables ${\bf{a}}_1, {\bf{a}}_2, \cdots$,
and ${\bf{a}}_m$ are perpendicular to each other.

Proof:

$r_{{\bf{a}}_1 {\bf{a}}_2 \cdots {\bf{a}}_m }^2=0 \Leftrightarrow
\omega_{{\bf{a}}_1 {\bf{a}}_2 \cdots {\bf{a}}_m }^2=1$

According to Properties 3.2 and 3.3, we can obtain
\[
1 = \omega _{{\bf{a}}_1 {\bf{a}}_2  \cdots {\bf{a}}_m }^2  \le
\cdots \le \omega _{{\bf{a}}_i {\bf{a}}_j }^2 \le 1
\]

Hence, for arbitrary $i,j\in {1,2,\cdots,m},i\neq j$, we have
$\omega _{{\bf{a}}_i {\bf{a}}_j }^2  = 1$ and $r _{{\bf{a}}_i
{\bf{a}}_j }^2 =\rho _{{\bf{a}}_i
{\bf{a}}_j }^2 = 0$.

\hfill$\blacksquare$

{\bf{Property 3.6}} \ If variables ${\bf{a}}_1, {\bf{a}}_2, \cdots,$
and ${\bf{a}}_{m-1}$ are not linear dependent,
$r_{{\bf{a}}_1{\bf{a}}_2\cdots{\bf{a}}_m}^2$ then gets the biggest
value 1 if and only if variable ${\bf{a}}_m$ lies on the hyperplane
spanned by ${\bf{a}}_1, {\bf{a}}_2, \cdots,$ and ${\bf{a}}_{m-1}$,
and $r_{{\bf{a}}_1{\bf{a}}_2\cdots{\bf{a}}_m}^2$ gets the smallest
value $r_{{\bf{a}}_1{\bf{a}}_2\cdots{\bf{a}}_{m-1}}^2$ if and only
if ${\bf{a}}_m$ is perpendicular to the hyperplane spanned by
${\bf{a}}_1, {\bf{a}}_2, \cdots,$ and ${\bf{a}}_{m-1}$.

Proof: The first half is true according to Property 3.4. Now we
prove the second part.

According to Lemma 2, we have
\[
 \gamma _{{\bf{a}}_1 {\bf{a}}_2  \cdots {\bf{a}}_m }^2  - \gamma _{{\bf{a}}_1 {\bf{a}}_2  \cdots {\bf{a}}_{m-1}
 }^2
 = \sum\limits_{J_{m - 2} } {(\sum\limits_{p \notin J_{m - 2} } \begin{array}{l}
 \!\!( - 1)^{g(p:J_{m - 2} )} a'_{m,p} \det [{\bf{a}}'(m - 1)|p,J_{m - 2} ] )^2 \\
 \end{array} }.
\]

Denote the determinant $\det (\mathop i\limits^ -  ;J_{m - 2} )$ as
\[
\det (\mathop i\limits^ -  ;J_{m - 2} ) = \det ([{\bf{a}}'_1
,{\bf{a}}'_2 , \cdots ,{\bf{a}}'_{m - 1} \backslash {\bf{a}}'_i
|J_{m - 2} ]).
\]

Then we have
\[
\begin{array}{l}
 \sum\limits_{p \notin J_{m - 2} } {( - 1)^{g(p:J_{m - 2} )} a'_{m,p} } \det [{\bf{a}}'(m - 1)|p,J_{m - 2} ] \\
 = \sum\limits_{p \notin J_{m - 2} } {a'_{m,p} } \sum\limits_i {( - 1)^{i + 1} a'_{i,p} \det (\mathop i\limits^ -  ;J_{m - 2}
 )} \\
  = \sum\limits_i {( - 1)^{i + 1} \det (\mathop i\limits^ -  ;J_{m - 2} )\sum\limits_{p \notin J_{m - 2}} {a'_{m,p} a'_{i,p} } }  \\
  = \sum\limits_i {( - 1)^{i + 1} \det (\mathop i\limits^ -  ;J_{m - 2} )(\sum\limits_p {a'_{m,p} a'_{i,p} }  - \sum\limits_{p \in J_{m - 2}} {a'_{m,p} a'_{i,p} } } ) \\
  = \sum\limits_i {( - 1)^{i + 1} \det (\mathop i\limits^ -  ;J_{m - 2} )({\bf{a}}'_i ,{\bf{a}}'_m
  )} - \sum\limits_i {( - 1)^{i + 1} \det (\mathop i\limits^ -  ;J_{m - 2} )} \sum\limits_{p = j_1 }^{j_{m - 2} } {a'_{m,p} a'_{i,p} }  \\
 \end{array}
\]

The first part of the above equation is a formula of the inner
product. The second part can be simplified as
\[
\begin{array}{l}
 \sum\limits_i {( - 1)^{i + 1} \det (\mathop i\limits^ -  ;j_1 ,j_2 , \cdots ,j_{m - 2} )} \sum\limits_{p = j_1 }^{j_{m - 2} } {a'_{m,p} a'_{i,p} }  \\
  = \det \!\!\left[ {\begin{array}{*{20}c}
   {\!\!\sum\limits_{p = j_1 }^{j_{m - 2} } {a'_{m,p} a'_{i,p} } } & \!\!\!{a'_{1,j_1 } } & \!\!\!{a'_{1,j_2 } } &  \!\!\!\cdots  & \!\!\!{a'_{1,j_{m - 2} } }  \\
   {\!\!\sum\limits_{p = j_1 }^{j_{m - 2} } {a'_{m,p} a'_{i,p} } } & \!\!\!{a'_{2,j_1 } } & \!\!\!{a'_{2,j_2 } } &  \!\!\!\cdots  & \!\!\!{a'_{2,j_{m - 2} } }  \\
    \!\!\vdots  &  \!\!\!\vdots  &  \!\!\!\vdots  &  \!\!\!\ddots  &  \!\!\!\vdots   \\
   {\!\!\sum\limits_{p = j_1 }^{j_{m - 2} } {a'_{m,p} a'_{i,p} } } & \!\!\!{a'_{m - 1,j_1 } } & \!\!\!{a'_{m - 1,j_2 } } &  \!\!\!\cdots  & \!\!\!{a'_{m - 1,j_{m - 2} } }  \\
\end{array}} \!\!\!\!\right]\!=0 \\
 \end{array}
\]

Hence,
\[
\begin{array}{l}
\sum\limits_{p \notin J_{m - 2} } {( - 1)^{g(p:J_{m - 2} )} a'_{m,p}
} \det [{\bf{a}}'(m - 1)|p,J_{m - 2} ] \\= \sum\limits_{i = 1}^{m -
1} {( - 1)^{i + 1} \det (\mathop i\limits^ -  ;J_{m - 2}
)({\bf{a}}'_i ,{\bf{a}}'_m )}.
\end{array}
\]

Then $\gamma _{{\bf{a}}_1 {\bf{a}}_2 \cdots {\bf{a}}_m }^2$ gets the
minimum value $r_{{\bf{a}}_1{\bf{a}}_2\cdots{\bf{a}}_{m-1}}^2$ if
and only if
\begin{equation}
\sum\limits_{i = 1}^{m - 1} {( - 1)^{i + 1} \det (\mathop i\limits^
-  ;j_1 ,j_2 , \cdots ,j_{m - 2} )({\bf{a}}'_i ,{\bf{a}}'_m )} =0
\end{equation}
holds for all possible $\{j_1 ,j_2 , \cdots ,j_{m - 2}\}$ and $i$.

Because variables ${\bf{a}}_1, {\bf{a}}_2, \cdots, {\bf{a}}_{m-1}$
are not linear dependent, $J_{m - 1}=\{j_1,j_2, \cdots$, $j_{m-1}$\}
exists to make the rank of $[{\bf{a}}'(m-1)|J_{m - 1}]$ equal to
$m-1$. Then $[{\bf{a}}'(m-1)|J_{m - 1}]$ is an invertible matrix and
its adjoint matrix is also an invertible matrix. Each row of the
adjoint matrix of $[{\bf{a}}'(m-1)|J_{m - 1}]$ is just the linear
coefficients of one equation in the above equation. Finally, we
obtain $({\bf{a}}'_i ,{\bf{a}}'_m )=0, i=1,2,\cdots,m-1$.\hfill$\blacksquare$

According to Lemma 1 and Definition 2, if the number of variables is
the same as their dimension, we have the following corollary, which
gives a new explanation of determinant from the view of multivariate
correlation.

\noindent{\bf{Corollary 3}} \ For $n$-dimensional non-zero-variance variables ${\bf{a}}_1, {\bf{a}}_2, \cdots$, and ${\bf{a}}_n$,
if ${\bf{a}}'_i$ is the standardized vector of ${\bf{a}}_i$,
$i\in\{1,2,\cdots,n\}$, and ${\bf{A}}=[{\bf{a}}'_1, {\bf{a}}'_2,
\cdots, {\bf{a}}'_n]$, then we have
\begin{equation}
\omega _{{\bf{a}}_1 {\bf{a}}_2  \cdots {\bf{a}}_m }^2 = (\det
({\bf{A}}))^2 .
\end{equation}

This corollary gives a new explanation of determinant that if the
row or column vectors of a matrix are all standardized, then the
absolute value of the determinant of the matrix depicts the linear
irrelevance of these standardized vectors.

Corollary 3 prompts us to consider the sign of the proposed UCC and
UIC. If UIC is defined by the determinant of ${\bf{A}}$ when the
number of variables is the same as their dimension, the value of
UIC will have a positive or negative sign for a group of variables.
However, because the sign of the determinant of ${\bf{A}}$ varies
with the order of these variable appeared in the matrix ${\bf{A}}$,
it is meaningless to take time to decide which sign is better.
Moreover, the absolute value of correlation coefficient is more appropriate to measure the
strength of correlation.

Lastly, if the variables in the inner product matrix ${\bf{M}}$ are
all standardized vectors, the inner product matrix is then transformed into
the correlation matrix. Correlation matrix is also a widely used
tool in various fields. For $n$-dimensional non-zero-variance
variables ${\bf{a}}_1$, ${\bf{a}}_2$, $\cdots$, ${\bf{a}}_m$, $2\le
m \le n$, if $\rho_{ij}$ is Pearson's correlation
coefficient between ${\bf{a}}_i$ and ${\bf{a}}_j$,
$i,j\in\{1,2,\cdots,m\}$, the correlation matrix ${\bf{R}}$ of these
variables is as following:
\[
{\bf{R}} = \left[ {\begin{array}{*{20}c}
   {\rho_{11}} & {\rho_{12}} &  \cdots  & {\rho_{1m } }  \\
   {\rho_{21 } } & {\rho_{22}} &  \cdots  & {\rho_{2m } }  \\
    \vdots  &  \vdots  &  \ddots  &  \vdots   \\
   {\rho_{m1 } } & {\rho_{m2 } } &  \cdots  & {\rho_{mm } }  \\
\end{array}} \right],
\]
in which the diagonal elements are all 1.

Then we have the following Multivariate correlation Theorem
from Lemma 1, Definition 1 and Definition 2:

\noindent{\bf{Multivariate Correlation Theorem}} \ If
$r_{{\bf{a}}_1{\bf{a}}_2 \cdots {\bf{a}}_m }$ and
$\omega_{{\bf{a}}_1{\bf{a}}_2 \cdots {\bf{a}}_m }$ are the unsigned correlation coefficient and the unsigned incorrelation coefficient for $n$-dimensional
non-zero-variance variables ${\bf{a}}_1$, ${\bf{a}}_2$,
$\cdots$, and ${\bf{a}}_m$, respectively, then
\begin{equation}
\begin{array}{l}
\omega_{{\bf{a}}_1{\bf{a}}_2 \cdots {\bf{a}}_m }^2 = \det
({\bf{R}}) \\
r_{{\bf{a}}_1{\bf{a}}_2 \cdots {\bf{a}}_m}^2 = 1 - \det
({\bf{R}})
\end{array}.
\end{equation}

Then according to the inner product-determinant equation, we have the UIC equation as following:

\noindent{\bf{UIC Equation}} \ For
$n$-dimensional non-zero-variance variables ${\bf{a}}_1$, ${\bf{a}}_2$, $\cdots$, and
${\bf{a}}_m$, $m\leq n$, 
\begin{equation}
\omega_{{\bf{a}}_1{\bf{a}}_2 \cdots {\bf{a}}_m }^2=\sum\limits_\pi  {2^{\left| {\pi _3 } \right|} ( - 1)^{m - \left|
{\pi} \right|} }{\prod\limits_{\pi _2 } \rho^2_{{\bf{a}}_i{\bf{a}}_j } }
\prod\limits_{\pi _3 } {\Re\!\!<\!{\bf{a}}_{k_1} ,{\bf{a}}_{k_2} ,
\cdots ,{\bf{a}}_{k_p }\!>}
\end{equation}
where $\Re\!\!<\!{\bf{a}}_{k_1} ,{\bf{a}}_{k_2} ,
\cdots ,{\bf{a}}_{k_p }\!>= \rho_{{\bf{a}}_{k_1}{\bf{a}}_{k_2}}   \rho_{{\bf{a}}_{k_2}{\bf{a}}_{k_3}}
\cdots \rho_{{\bf{a}}_{k_{p-1}} {\bf{a}}_{k_p}}\rho_{{\bf{a}}_{k_p}{\bf{a}}_{k_1}}$ is the circular correlation coefficient, and $\pi$, $\pi_1$, $\pi_2$, $\pi_3$, and $\left|
{\pi _3 } \right|$ are kept the same meanings as that in the inner product-determinant equation in which the inner products are replaced by correlation coefficients. 

\vspace{0.2 cm}

{\center\bf{4. \quad VISUALIZATION ANALYSIS}

}

Here we take the unsigned tri-variate correlation
coefficient as an example to visually display the proposed UCC and discuss its effectiveness.

\vspace{0.1cm} \noindent{\bf{4.1 \ Visualization of the Proposed
UCC}}\vspace{0.1cm}

\begin{figure*}[t]
\vspace{2mm} \centering \centerline{\epsfig{figure=Fig3.pdf,width=
5.9 in}} \caption{Some surfaces of the unsigned correlation
coefficient against the angles $\beta$ and $\gamma$ with the other
angle $\alpha$ fixed as $30^0$, $90^0$, $145^0$, and $160^0$,
respectively.} \label{Fig3}
\end{figure*}

\begin{figure*}[t]
\centering \centerline{\epsfig{figure=Fig4.pdf,width= 12.5 cm}}
\caption{Some equipotential lines with different unsigned correlation coefficients
for three variables with fixed the angle $\alpha = 90^0$ and $\alpha
= 140^0$, respectively.} \label{Fig4}
\end{figure*}

As shown in Fig. 1, we can easily show the effectiveness of Pearson's correlation coefficient by a 2D graph. For UCC with more variables, it is impossible to visually display the relation among UCC and these spatial angles in a 2D or 3D graph. However, if one angle $\alpha$ is fixed, then the relation of the other two
angles $\beta$ and $\gamma$, and the value of UCC among three variables can be visually display in a 3D graph. Four such
graphs are shown in Figure~\ref{Fig3} with different fixed angles
$\alpha=30^0$, $\alpha=90^0$, $\alpha=145^0$, and $\alpha=160^0$,
respectively. From Figure~\ref{Fig3} we can see that these surfaces
have the similar structure but different depth, curvature, and top
rectangles.

\vspace{0.1cm} \noindent{\bf{4.2 \ Contour Line and Geometrical
Explanation of UIC}}\vspace{0.1cm}

From Figure~\ref{Fig3} we can see that a myriad of contour lines exist in
the surfaces. A simple example is that if the
pairwise angles for three variables ${\bf{a}}_1$, ${\bf{b}}_1$, and
${\bf{c}}_1$ are $45^0$, $45^0$, and $60^0$, respectively, and the
pairwise angles for three variables ${\bf{a}}_2$, ${\bf{b}}_2$, and
${\bf{c}}_2$ are $30^0$, $90^0$, and $90^0$, respectively, then the
three variables ${\bf{a}}_1$, ${\bf{b}}_1$, and ${\bf{c}}_1$ and another three variables ${\bf{a}}_2$, ${\bf{b}}_2$, and ${\bf{c}}_2$
have the same value of UCC. Some equipotential lines for three variables
with the fixed angle $\alpha$ equal to $90^0$ and $140^0$ are shown
in Figure~\ref{Fig4}.

In fact, the linear relation for multiple variables is closely tied
with the parallelotope in multi-dimensional space. For example, the
linear space structured by $m$ independent variables is the
$m$-dimensional linear space, and the vector sum of the $m$
variables is the diagonal of the parallelogram formed by these
variables in this $m$-dimensional space.

Barth had proposed that the determinant of a Gram matrix is the
square of the volume of the parallelotope formed by the variables
\cite{Barth}, and the correlation matrix is a special Gram matrix.
According to Section 3 in this paper, the square of the proposed
unsigned incorrelation coefficient is the determinant of correlation
matrix. Hence, we have the following corollarys:

\noindent{\bf{Corollary 4}} \ If $\omega_{{\bf{a}}_1 {\bf{a}}_2 \cdots
{\bf{a}}_m }$ is the unsigned incorrelation coefficient among multiple
variables ${\bf{a}}_1$, ${\bf{a}}_2$, $\cdots$, ${\bf{a}}_m$, then
$\omega_{{\bf{a}}_1 {\bf{a}}_2 \cdots {\bf{a}}_m }$ is the volume of
the parallelotope formed by the vectors ${\bf{a}}'_1$,
${\bf{a}}'_2$, $\cdots$, ${\bf{a}}'_m$, where ${\bf{a}}'_i$ is the
standardized vector of ${\bf{a}}_i, i=1,2,\cdots,m$.

\noindent{\bf{Corollary 5}} \ For spatial angles $\alpha_1, \alpha_2,
\cdots, \alpha_{\text{C}^2_m}$, the unsigned incorrelation coefficient is the
volume of the parallelotope formed by the unit vectors whose
pairwise angles are $\alpha_1, \alpha_2$, $\cdots$, and $\alpha_{\text{C}^2_m}$,
respectively.

According to Corollary 5, the equipotential lines on the UCC surface is
the equal-volume line for different spatial angles.

\vspace{0.1cm} \noindent{\bf{4.3 \ Effectiveness of the Proposed
UCC}}\vspace{0.1cm}

\begin{figure}[t]
\centering \centerline{\epsfig{figure=Fig5.pdf,width= 7.3 cm}}
\caption{A case of $\alpha,\beta$, and $\gamma$. $\alpha, \beta$, and $\gamma$ are the angles between ${\bf{a}}$ and
${\bf{b}}$, between ${\bf{b}}$ and ${\bf{c}}$, and between
${\bf{a}}$ and ${\bf{c}}$, respectively. $\left\|{\text{OA}}\right\|=\left\|{\text{OB}}\right\|=\left\|{\text{OC}}\right\|=1$.
Points C$_0$ is the projection point of C on the plane spanned by ${\bf{a}}$ and ${\bf{b}}$. $\eta=\angle$COC$_0$. }
\label{Fig5}
\end{figure}

In Section 3, some important properties have shown that the proposed UCC and UIC are effective measures for correlation of multivariate variables. Here we visually verify the effectiveness of UCC for three variables.

The effective of Pearson's correlation coefficient has been visually displayed in Figure 1, from which we can see that the strength of correlation measured by correlation coefficient is very close to the strength of correlation measured by the angle $\gamma$. In subsection 4.1, several 3D figures of the unsigned tri-variate correlation
coefficient with a fixed angle $\alpha$ are provided, the effectiveness of which will be discussed here.

In fact, we can take this case of the unsigned tri-variate correlation
coefficient with a fixed $\alpha$ 
as that ${\bf{a}}$ and ${\bf{b}}$ are fixed and ${\bf{c}}$ can be any vectors starting from O. Then the strength of correlation among ${\bf{a}}$, ${\bf{b}}$, and ${\bf{c}}$ can be measured by the angle $\eta$ between ${\bf{c}}$ and OC$_0$, in which C$_0$ is the projection point of C on the plane spanned by ${\bf{a}}$ and ${\bf{b}}$. %the plane spanned by ${\bf{a}}$ and ${\bf{b}}$. %As discussed in the introduction part of this paper, the absolute value of Pearson's correlation coefficient can effectively measure bivariate correlation because it is always kept the opposite trend of the change of $\gamma^*$. Here we discuss whether the value of the unsigned tri-variate correlation coefficient is also a good measure of correlation when $\alpha$ is fixed.

We select three points A, B, and C on ${\bf{a}},{\bf{b}}$, and ${\bf{c}}$, respectively, to make $\left\|{\text{OA}}\right\|=\left\|{\text{OB}}\right\|=\left\|{\text{OC}}\right\|=1$. Let UCC and UIC among ${\bf{a}},{\bf{b}}$, and ${\bf{c}}$ be $r_{{\bf{abc}}}$ and $\omega_{{\bf{abc}}}$, respectively. According to Corollary 4 and Corollary 5, the volume of the parallelotope formed by $\left\|{\text{OA}}\right\|, \left\|{\text{OB}}\right\|, \left\|{\text{OC}}\right\|$, and their other nine parallel line segments is $\omega_{{\bf{abc}}}$. Hence, we have
\begin{equation}
\sin\alpha\sin\eta=\omega_{{\bf{abc}}},
\end{equation}
then
\begin{equation}
r_{{\bf{abc}}} = \sqrt{1-\sin^2\alpha\sin^2\eta}\ .
%\cos \eta  = \sqrt {1 - \frac{{1 - {r^2}}}{1-{{{\cos }^2}\alpha }}}\ .
\end{equation}

Then the curves of the unsigned tri-variate correlation
coefficient against the angle $\eta$ are depicted in Figure~\ref{Fig6}, from which we can see that the strength of correlation among ${\bf{a}},{\bf{b}}$, and ${\bf{c}}$ measured by the unsigned tri-variate correlation coefficient (the black curves in Figure~\ref{Fig6}) is very close to the strength of correlation measured by the angle $\eta$ (the blue dotted segments in Figure~\ref{Fig6}).

\begin{figure}[t]
\centering \centerline{\epsfig{figure=Fig6.pdf,width= 15.1 cm}}
\caption{The curves of the unsigned tri-variate correlation
coefficient against the angle $\eta$. The angle $\alpha$ is fixed with different values, and $\eta$ is $\angle$COC$_0$ in Figure~\ref{Fig5}.}
\label{Fig6}
\end{figure}

%Therefore, $r_{{\bf{abc}}}$ is a continuous and strictly monotone decreasing function of $\eta$.

\vspace{0.2 cm}

{\center\bf{5. \quad CONCLUSION}

}

\end{document}